%Version 1.1, 9/3/03, omission of ``fiber functor'' in introduction corrected
%1.01, 8/10/03
%Version 1.0, 8/5/03 (first arXiv version)

%Determining a semisimple group from its representation degrees

%pseudo-Israeli macros
\catcode`\@=11
\def\undefine#1{\let#1\undefined}
\def\newsymbol#1#2#3#4#5{\let\next@\relax
 \ifnum#2=\@ne\let\next@\msafam@\else
 \ifnum#2=\tw@\let\next@\msbfam@\fi\fi
 \mathchardef#1="#3\next@#4#5}
\def\mathhexbox@#1#2#3{\relax
 \ifmmode\mathpalette{}{\m@th\mathchar"#1#2#3}%
 \else\leavevmode\hbox{$\m@th\mathchar"#1#2#3$}\fi}
\def\hexnumber@#1{\ifcase#1 0\or 1\or 2\or 3\or 4\or 5\or 6\or 7\or 8\or
 9\or A\or B\or C\or D\or E\or F\fi}

\newdimen\ex@
\ex@.2326ex
\def\varinjlim{\mathop{\vtop{\ialign{##\crcr
 \hfil\rm lim\hfil\crcr\noalign{\nointerlineskip}\rightarrowfill\crcr
 \noalign{\nointerlineskip\kern-\ex@}\crcr}}}}
\def\varprojlim{\mathop{\vtop{\ialign{##\crcr
 \hfil\rm lim\hfil\crcr\noalign{\nointerlineskip}\leftarrowfill\crcr
 \noalign{\nointerlineskip\kern-\ex@}\crcr}}}}
\def\varliminf{\mathop{\underline{\vrule height\z@ depth.2exwidth\z@
 \hbox{\rm lim}}}}

\font\tenmsa=msam10
\font\sevenmsa=msam7
\font\fivemsa=msam5
\newfam\msafam
\textfont\msafam=\tenmsa
\scriptfont\msafam=\sevenmsa
\scriptscriptfont\msafam=\fivemsa
\edef\msafam@{\hexnumber@\msafam}
\mathchardef\dabar@"0\msafam@39
\def\dashrightarrow{\mathrel{\dabar@\dabar@\mathchar"0\msafam@4B}}
\def\dashleftarrow{\mathrel{\mathchar"0\msafam@4C\dabar@\dabar@}}

\font\tenmsb=msbm10
\font\sevenmsb=msbm7
\font\fivemsb=msbm5
\newfam\msbfam
\textfont\msbfam=\tenmsb
\scriptfont\msbfam=\sevenmsb
\scriptscriptfont\msbfam=\fivemsb
\edef\msbfam@{\hexnumber@\msbfam}
\def\Bbb#1{{\fam\msbfam\relax#1}}
\def\widehat#1{\setbox\z@\hbox{$\m@th#1$}%
 \ifdim\wd\z@>\tw@ em\mathaccent"0\msbfam@5B{#1}%
 \else\mathaccent"0362{#1}\fi}
\font\teneufm=eufm10
\font\seveneufm=eufm7
\font\fiveeufm=eufm5
\newfam\eufmfam
\textfont\eufmfam=\teneufm
\scriptfont\eufmfam=\seveneufm
\scriptscriptfont\eufmfam=\fiveeufm

\newsymbol\boxtimes 1202
\newsymbol\ndiv 232D
\newsymbol\supsetneqq 2325
\newsymbol\subsetneqq 2324
\newsymbol\ngeq 2303

\catcode`\@=12

%My fonts
\magnification=\magstep1
\font\title = cmr10 scaled \magstep2

\font\smalltext = cmr7
\font\smallmath= cmmi7
\font\tinymath=cmmi5
\font\smallsym = cmsy7
\font\author = cmcsc10
\font\addr = cmti7
\font\byabs = cmr7

%format information
\parindent=1em
\baselineskip 15pt
\hsize=12.3 cm
\vsize=18.5 cm

%My format macros
\newcount\refcount
\newcount\seccount
\newcount\sscount
\newcount\eqcount
\newcount\boxcount
\newcount\testcount
\newcount\bibcount
\boxcount = 128
\seccount = -1

\seccount = -1

\def\sec#1{\advance\seccount by 1\bigskip\goodbreak\noindent
	{\bf\number\seccount.\ #1}\medskip \sscount = 0\eqcount = 0}
\def\subsec{\advance\sscount by 1\medskip\goodbreak\noindent
	{{\number\seccount.\number\sscount}\ \ } \eqcount = 0}
\def\namess#1{\advance\sscount by 1\global
	\edef#1{\number\seccount.\number\sscount}
	\medskip\goodbreak\noindent{\bf #1\ \ } \eqcount = 0}
\def\proc#1#2{\advance\sscount by 1\eqcount = 0
	\medskip\goodbreak\noindent{\author #1}
	{\tenrm{\number\seccount.\number\sscount}}:\ \ {\it #2}}
\def\nproc#1#2#3{\advance\sscount by 1\eqcount = 0\global
	\edef#1{#2\ \number\seccount.\number\sscount}	
	\medskip\goodbreak\noindent{\author #2}
	{\tenrm{\number\seccount.\number\sscount}}:\ \ {\it #3}}
\def\proof{\medskip\noindent{\it Proof:\ \ }}

\def\eql#1{\global\advance\eqcount by 1\global
	\edef#1{(\number\seccount.\number\sscount.\number\eqcount)}\leqno{#1}}
\def\ref#1#2{\advance\refcount by 1\global
	\edef#1{[\number\refcount]}\setbox\boxcount=
	\vbox{\item{[\number\refcount]}#2}\advance\boxcount by 1}
\def\biblio{{\frenchspacing
	\bigskip\goodbreak\centerline{\bf REFERENCES}\medskip
	\bibcount = 128\loop\ifnum\testcount < \refcount
	\goodbreak\advance\testcount by 1\box\bibcount
	\advance\bibcount by 1\vskip 4pt\repeat\medskip}}

%Israel J. Macros
\def\emph#1{{\bf #1}}
\def\colon{{:}\;}
\def\|{|\;}

%My math symbols

\def\scirc{{\scriptstyle\circ}}

\def\C{{\Bbb C}}
\def\R{{\Bbb R}}
\def\N{{\Bbb N}}
\def\Z{{\Bbb Z}}
\def\Q{{\Bbb Q}}
\def\F{{\Bbb F}}

\def\rt#1{{\buildrel #1\over\longrightarrow}}

\def\Aut{{\rm Aut}}
\def\Out{{\rm Out}}

\def\GL{{\rm GL}}

\def\qed{\hfill\hbox{$\sqcup$\llap{$\sqcap$}}\medskip}

\def\PSp{{\rm PSp}}
\def\SO{{\rm SO}}
\def\SU{{\rm SU}}

\def\tr{{\rm tr}}

\def\Span{{\rm Span}\,}

\def\eff{{\rm eff}}
\def\lev{{\rm lev}}
\def\ord{{\rm ord}}
\def\mod#1{(\mathop{\rm mod}\nolimits #1)}

\def\Rep{{\rm Rep}}
\def\hf{{1\over 2}}

\ref\Bourbaki{Bourbaki, N.: Groupes et alg\`ebres de Lie, Chap. 5--6 (1968),
	Chap. 7--8 (1975), Hermann, Paris.}

\ref\Tannaka{Deligne, Pierre; Milne, James S.; Ogus, Arthur; Shih, Kuang-yen: Hodge cycles,
	motives, and Shimura varieties. Lecture Notes in Mathematics, 900.  Springer-Verlag, 
	Berlin-New York,  1982.}

\ref\Gassmann{Gassmann, F.: Bemerkungen zu vorstehenden Arbeit von Hurwitz.
	{\it Math. Z.} {\bf 25} (1926) 124--143.}

\ref\simpLie{Moody, R.~V.; Patera, J.; Rand, D.~W.:
	simpLie Version 2.1, Macintosh Software for Representations of Simple Lie Algebras.}

\ref\Pink{Larsen, M.; Pink, R.: Determining representations from invariant dimensions.  
	{\it Invent. Math.}  {\bf 102}  (1990),  no. 2, 377--398.}

\ref\Sunada{Sunada, T.: Riemannian coverings and isospectral manifolds.
	{\it Annals of Math.} {\bf 121} (1985) 169--186.}

\ref\Witten{Witten, E.: On quantum gauge theories in two dimensions. 
	{\it Comm. Math. Phys.} {\bf 141}  (1991),  no. 1, 153--209.}

\centerline{\title Determining a Semisimple Group}
\medskip
\centerline{\title from its Representation Degrees}
\bigskip
\centerline{\byabs BY}
\medskip
\noindent{\author\hfill Michael Larsen\footnote*
{\tenrm Partially supported by NSF
Grant DMS-0100537}\hfill}
\medskip
\centerline{\addr Department of Mathematics, Indiana University}
\centerline{\addr Bloomington, IN 47405, USA}
\bigskip

\centerline{\byabs ABSTRACT}
\smallskip
{\byabs \narrower\narrower
\textfont0 = \smalltext
\textfont1 = \smallmath
\scriptfont1 = \tinymath
\textfont2 = \smallsym
The Lie algebra of a compact semisimple Lie group $G$ is determined by the degrees of the irreducible representations of $G$.  However, two different groups can have the same representation degrees.
\par}

\medskip
%-------------------------------------------------------------------------
\sec{Introduction}

A compact semisimple Lie group $G$ has a finite number of representations of any given degree.  In this paper we investigate to what extent $G$ is determined by the multi-set of degrees of its representations or equivalently, the multi-set of degrees of its irreducible representations.
This problem arises naturally in at least two different ways.  On the one hand, following
Witten \Witten, we define the zeta-function $\zeta_G(s)$ of a compact semisimple
Lie group as the Dirichlet series
$$\zeta_G(s)=\sum_{V\;{\rm irreducible}}(\dim V)^{-s}.$$
The question in this paper then belongs to the familiar class of problems asking when two different mathematical objects can have the same zeta-function.  On the other hand, the multi-set of degrees of the irreducible representations of $G$ comprises a (small) part of the data encoded in the representation category $\Rep_{\C}G$.  We know \Tannaka\ that the full structure of neutral Tannakian category on $\Rep_{\C}G$ enables us to recover the complexification $G_{\C}$ of $G$ and therefore $G$ itself (as a maximal compact subgroup of $G_{\C}$).   For semisimple groups, though, there are results (see for instance, \Pink) showing that $\Rep_{\C}G$ has more data than we need to determine $G$.  In this paper we are asking whether the $\C$-linear abelian category 
with fiber functor underlying $\Rep_{\C}G$ already determines $G$.

The main theorem of this paper is that $\zeta_G(s)$ determines the Lie algebra of $G$ up to isomorphism.   Our basic tool is the Weyl dimension formula, which expresses the dimension of each irreducible representation of a group $G$ as a product over the set of positive roots of $G$.   We try to recover the geometry of the root system of $G$ from the factorizations of the representation degrees, the basic difficulty being that there are many ways to factor a given degree.  Our strategy is to choose values of $n$ which have few factorizations or at least few factorizations which could possibly arise from Weyl's formula.

The simplest idea is to consider prime powers $n = p^k$.   It is an easy consequence of Weyl's formula that $\dim V_{p\lambda+(p-1)\rho} = p^{|R^+|} \dim V_\lambda$, where $V_\lambda$ is the irreducible representation of $G$ of highest weight $\lambda$, $R^+$ is the set of positive roots of $G$, and $\rho$ is the half-sum of the elements of $R^+$.   An elaboration of this argument, given in \S1,
shows that the prime power coefficients determine the number of irreducible factors of $R$ and the number of roots in each factor.  Unfortunately, there are infinitely many pairs of distinct irreducible root systems with the same number of roots, and prime power coefficients alone are in general insufficient to disambiguate further.  For example, $\SO(15)$ and $\PSp(14)$ both have the property that the number of their irreducible representations of degree $p^k$ is $1$ if $k\in 49\Z$ and $0$ otherwise (for $p\gg 0$, this follows from the discussion in \S1, but a computer search is needed to confirm it for small primes.) 

We therefore broaden our search to encompass degrees $n$ in which the largest power of some prime $p$ dividing $n$ is not too much smaller than $n$.  To illustrate this idea, we consider the case that $R$ is irreducible, and exclude weights $\lambda$ with $\lambda+\rho$ divisible by $p$ (such weights are not ``allowable at $p$'' in the terminology of this paper).  We compare the size of the logarithm of the largest power of $p$ dividing $\dim V_\lambda$ to that of $\dim V_\lambda$ itself
and prove that
$$\lim_{p\to\infty}\sup_{\lambda\neq 0\;\rm{allowable}}{-\log|\dim V_\lambda|_p\over \log\dim V_\lambda}
= \lim_{p\to\infty}\sup_{\lambda\neq 0\;\rm{allowable}}{\log p\;\ord_p\dim V_\lambda\over \log\dim V_\lambda}\eql\limitDef$$
approaches a rational limit strictly between $0$ and $1$ which depends on $R$.  This limit is the \emph{efficiency} of $R$ (defined in \S3 directly in terms of the geometry of $R$).  A more efficient root system always has allowable representation degrees which cannot be achieved by an allowable weight of a less efficient system.   Once the efficiency is known, we can begin studying the families of representations (indexed by $p$) whose degrees achieve it.    It turns out the optimal ones, in a sense made explicit below, lie in finitely many one parameter families of the form $p\mapsto p\mu+\nu$.
For each such family, the degree is given as a polynomial in $p$, the ``Weyl polynomial'' associated to $\mu$ and $\nu$.  The Weyl polynomials of a simple root system encode enough information, ultimately, to extract the root system giving rise to them.  

The process of excluding non-allowable representations is carried out in \S1 and is somewhat delicate; because our groups may not be products of almost simple factors, we must restrict attention to primes which are congruent to $1$ modulo a sufficiently divisible integer.   Weyl polynomials are introduced in \S2 and the basic pairs $(\mu,\nu)$ used in the proof of the main theorem are given.   The computation of efficiencies occupies much of \S3 and is rather involved.  There does not seem to be a simple formula for the answer: for instance, the efficiency of $E_8$ turns out to be $63\over 117$.   
The proof of the main theorem, by means of a rather complicated induction on the number of factors, is given in \S4.

To show that our result is in some sense sharp, in \S5 we give a construction of pairs of non-isomorphic compact semisimple Lie groups with the same Witten zeta-function.
This construction is modelled on a well known theorem of F.~Gassmann \Gassmann\ asserting that the Dedekind zeta-function $\zeta_K(s)$ of a number field $K$
does not determine $K$ up to isomorphism.   The proof of Gassmann's theorem is purely group-theoretic.  Namely, there exists a finite group $G$ with subgroups $H_1$ and $H_2$, not conjugate in $G$, whose elements can be put into
one-to-one correspondence in such a way that corresponding elements {\it are} conjugate in $G$.
The fixed fields of $H_1$ and $H_2$ in a Galois extension with group $G$ then have the same 
zeta-function.  T.~Sunada \Sunada\ used the same trick to construct pairs of 
non-isometric isospectral manifolds; here $\pi_1$ plays the role of the Galois group.
Like Sunada, we exploit $\pi_1$,
but our fundamental  groups are abelian, so there are no inner automorphisms to work with.
Instead we make use of outer automorphisms.
We construct subgroups $H_1$ and $H_2$ of a suitable $\pi_1(G)$, whose elements can be put into one-to-one correspondence
in such a way that corresponding elements always lie in the same orbit of $\Aut(\pi_1(G))$
but such that the groups $H_1$ and $H_2$ as a whole do not.
The corresponding covering groups can then be shown to have the same representation degrees.

This work was initiated during a visit to the Hebrew University in Jerusalem in the summer of 2002.  I would like to thank the Einstein Institute of Mathematics for its hospitality during this visit.
\medskip

%---------------------------------------------------------------
\sec{Modifying the zeta-function}

Let $G$ be a compact semisimple group, $\tilde G$ its universal cover, and $G_1,\ldots, G_m$
the almost simple direct factors of $G$.  Let $\Lambda_G$ and $\Lambda_{\tilde G}$ denote
the weight lattices of $G$ and $\tilde G$ respectively, so 
$$[\Lambda_{\tilde G}:\Lambda_G] = |\ker(\tilde G\to G)|.$$
We write $V_\lambda$ for the irreducible representation of $\tilde G$ with highest weight $\lambda$;
we can regard $V_\lambda$ as a $G$-representation if and only if $\lambda\in\Lambda_G$,
and all irreducible representations of $G$ arise in this way.  We have a direct sum decomposition
$$\Lambda_{\tilde G}=\bigoplus_{i=1}^m \Lambda_{G_i},$$
and we write $\lambda=\lambda_1+\ldots+\lambda_m$.  Thus $\lambda$ 
is a dominant weight of $\tilde G$ if and only if $\lambda_i$ is a dominant weight of $G_i$
for all $i$.   There is a tensor decomposition
$$V_\lambda = V_{\lambda_1}\otimes\cdots\otimes V_{\lambda_m}.$$
Let $R$ (resp. $R_i$) denote the root system of $G$ (resp. $G_i$).
We fix a positive Weyl chamber of $R$ (equivalently, positive Weyl chambers for each $R_i$.)
For each root $\alpha$, we denote the dual root $\alpha^\vee = {2\alpha\over \Vert\alpha\Vert^2}$, and we write $\rho$ (resp. $\rho_i$)
for the half-sum of positive roots in $R$ (resp. $R_i$).  The Weyl dimension formula
\Bourbaki~VIII~\S9~Th. 2 asserts
$$\dim V_\lambda = \prod_{\alpha\in R^+}{\alpha^\vee(\lambda+\rho)\over\alpha^\vee(\rho)}
= \prod_{i=1}^m\prod_{\alpha\in R_i^+}{\alpha^\vee(\lambda_i+\rho_i)\over\alpha^\vee(\rho_i)}
= \prod_{i=1}^m \dim V_{\lambda_i}.$$
It follows that 
$$\dim V_{p\lambda_i+(p-1)\rho_i}=p^{|R_i^+|}\dim V_{\lambda_i}.\eql\scaling$$
Thus,
$$\dim V_{(p^{a_i}-1)\rho_1+\cdots+(p^{a_m}-1)\rho_m} = p^{\sum a_i|R_i^+|}.$$
Note that when $p$ is odd, $(p^{a_i}-1)\rho_i$ lies in the root lattice of $G_i$, so
$$\sum_{i=1}^m(p^{a_i}-1)\rho_i \in \Lambda_G.$$
Therefore, every coefficient of $\zeta_G(s)$ is greater than equal to the corresponding coefficient of
$$\prod_{i=1}^m\bigl(1-p^{-|R_i^+|s}\bigr)^{-1}.$$
We will see that for fixed $G$, when $p\gg 0$ is a sufficiently large prime, the
latter expression gives precisely the part of the Dirichlet series $\zeta_G(s)$
consisting of $p$-power terms.

\nproc\Spanning
{Lemma}{An irreducible root system $R$ cannot be contained in the union of two proper subspaces of $\Span R$.}

\proof
Let $E = \Span R$ and $V$ and $W$ proper subspaces such that $R\subset V\cup W$.
We use induction on $\dim E$, the lemma being trivial when $\dim E = 1$.
Replacing $V$ and $W$ by $\Span (V\cap R)$ and $\Span (W\cap R)$ respectively, we may assume without loss of generality that $V\cap R$ are $W\cap R$
are root systems with span $V$ and $W$ respectively.
If $W\subset V$, then $R\subset V\subsetneqq E = \Span R$, which is absurd; likewise, $V$ is not a subspace of $W$.   If $W\subset V^\perp$,
then $R = (R\cap V) \cup (R\cap V^\perp)$, contrary to the hypothesis of irreducibility.
By the induction hypothesis, there exists a root $\alpha$ in $W$ which is neither in $V$ nor in $V^\perp$.
Thus $\alpha^\perp\cap V$ is a proper subspace of $V$, and there exists a root $\beta\in V\cap R$
which is neither in $W$ nor in $\alpha^\perp$.   As $\alpha$ and $\beta$ are neither scalar multiples of one another nor mutually orthogonal, either $\alpha+\beta$ or $\alpha-\beta$ belongs to $R$ (depending on whether the inner product of $\alpha$ and $\beta$ is positive or negative.)
Either way, $\alpha+\beta\not\in V\cup W$, and the lemma follows by induction.
\qed

\proc{Lemma}{If $p\gg 0$, then
$$\dim V_{\lambda_i}\in p^{\N}\Leftrightarrow\lambda_i+\rho_i\in p^{\N}\rho_i.$$
}

\proof
One direction is trivial.  By \scaling, it suffices to show $\dim V_{\lambda_i}\in p^{\N}$
implies $\lambda_i=0$ or $\lambda_i+\rho_i\in p\Lambda_{G_i}$.  If $p$ does not divide
$\lambda_i+\rho_i$, then
$$S=\{\alpha^\vee\in R_i^\vee\mid \alpha^\vee(\lambda_i+\rho_i)\in p\Z\}$$
is a proper root subsystem of $R_i^\vee$.   For $p\gg 0$, every root in $\alpha^\vee\in R\cap\Span S$ such that $p\alpha^\vee\in S$ lies in $S$ itself.  (In fact, an examination of all equal-rank subsystems of all irreducible root systems show that it suffices to take $p\ge 5$.)  Therefore, $S$ is closed.
By \Spanning,
the complement of $S$ spans $R_i^\vee$, so there are finitely many elements $\mu\in\Lambda_{G_i}$ 
such that
$$\beta^\vee(\mu+\rho_i)\le\prod_{\alpha^\vee\in R_i^\vee}\alpha^\vee(\rho_i)$$
for all $\beta^\vee\in R_i^\vee\setminus S$.  If $p\gg 0$, then $p\ndiv \dim V_\mu$ for all $\mu$ in this set.  Thus $\dim V_\mu\in p^{\N}$ implies $\mu=0$.
\qed

\nproc\Sizes
{Corollary}
{The zeta-function of $G$ determines, for each integer $n$, the number
of irreducible factors of the root system of $G$ of cardinality $n$.}
\medskip
This is not enough to determine the root system of $G$.
For one thing, the root systems $B_i$ and $C_i$ each have $2i^2$ roots.
We also have coincidences with $|R_i|=12$ ($A_3$ and $G_2$) $|R_i|=72$ ($B_6/C_6$, and $E_6$),
$|R_i|=240$ ($A_{15}$ and $E_8$), and a sparse but infinite set of cases in which a
root system of type A has the same number of roots of one of type B, C, or D.
So we need to refine the method.  We note, however, that fixing $\zeta_G(s)$ determines the
rank of $G$.  We may therefore define $N_G$ to be the factorial of the number of roots 
in the root system of $G$.  Thus $N_G$ depends on the one hand 
only on the Lie algebra of $G$ and on the other
only on the zeta-function of $G$.   It also has the property that $N_G\Lambda_G$ is contained in the root lattice of $G$.  We  define
$$\zeta_G^*(s)=\zeta_G(s)\prod_{p\equiv 1\ ({\rm mod}\;N_G)} \prod_{i=1}^m (1-p^{-|R_i^+|s}).$$
Note that for $p\equiv 1\pmod{N_G}$, 
$$\lambda_1+\cdots+\lambda_m\in\Lambda_G\Leftrightarrow\lambda_1+\cdots+\lambda_{i-1}
+(p\lambda_i+(p-1)\rho_i)+\lambda_{i+1}+\cdots+\lambda_m\in \Lambda_G.$$
\proc{Definition}
{We say that a weight $\lambda_i$ is \emph{allowable at a prime} $p$ if $p\ndiv \lambda_i+\rho_i$.}

\proc{Definition}
{We say that a weight $\lambda_1+\ldots+\lambda_m\in\Lambda_G$ is \emph{allowable}
if for every prime $p\equiv 1 \pmod{N_G}$ and every positive integer $i\le m$,
$\lambda_i$ is allowable at $p$.
We write $\Lambda_G^*$ for the set of allowable dominant weights of $G$.
%
%$$\Lambda_G^*=\{\lambda_1+\cdots+\lambda_m\in\Lambda_G\mid
%p\equiv 1\ \mod{N_G}\Rightarrow p\ndiv \lambda_i+\rho_i\ \forall i\},$$
%
Thus,
$$\zeta_G^*(s) = \sum_{\lambda\in\Lambda_G^*}\dim(V_\lambda)^{-s}.$$
}

It is easier to work with $\zeta_G^*(s)$ than with $\zeta_G(s)$ because for $\lambda\in\Lambda_G^*$,
the condition that $\dim V_\lambda$ is highly divisible by $p\equiv 1\ \mod{N_G}$ 
gives substantial information about the structure of
the Lie algebra of $G$, whereas for non-allowable $\lambda$, the same divisibility
condition could hold for essentially trivial reasons.  

%---------------------------------------------------------------
\sec{Weyl polynomials}

In this section we consider a fixed simply connected almost simple group $G$, i.e., $G$ is connected and simply connected, and $G/Z(G)$ is simple.
We are interested in the dimensions of representations of $G$ belonging to one-parameter 
families, specifically the families $p\mapsto p\mu+\nu$.
We let $R$ denote the root system of $G$; by the weight lattice of $R$, we mean the weight lattice of $G$.

\proc{Definition}{Given an irreducible root system $R$ and elements $\mu$ and $\nu$ of the weight lattice of $R$, the \emph{Weyl polynomial} $P^R_{\mu,\nu}(x)\in\Q[x]$ is given by the formula
$$P^R_{\mu,\nu}(x) = {\prod_{\alpha\in R^+}\bigl( \alpha^\vee(\mu)x + \alpha^\vee(\nu+\rho) \bigr)\over
\prod_{\alpha\in R^+}\alpha^\vee(\rho)}.$$
}

When $\lambda = n\mu  +\nu$ is dominant, 
$$\dim V_\lambda = P^R_{\mu,\nu}(n)\eql\eval$$
We will always assume that the following two conditions hold:
$$\alpha^\vee(\mu)\ge 0\ \forall\alpha\in R^+\eql\dominantI$$
and
$$\alpha^\vee(\mu) = 0\Rightarrow \alpha^\vee(\nu)\ge 0\ \forall\alpha\in R^+.\eql\dominantII$$
This implies that \eval\ holds for all $n\gg 0$.  If $p$ is a prime greater than 
the maximum of $|\alpha^\vee(\mu)|$, $|\alpha^\vee(\nu+\rho)|$, and $|\alpha^\vee(\rho)|$ for any root $\alpha$, then $\ord_p(\dim V_{p\mu + \nu})$ is equal to $|(\nu+\rho)^\perp\cap R^+|$, the multiplicity of $0$ as a root of $P^R_{\mu,\nu}(x)$.  In particular, this multiplicity is $|R^+|$ if and only
if $\nu+\rho = 0$, which corresponds to the case that for all $p$, $p\mu + \nu$ is not allowable at $p$.

\nproc\Explicit{Proposition}{For every irreducible root system $R$, there exist weights $\mu$ and $\nu$ satisfying \dominantI\ and \dominantII\ and such that
$(\nu+\rho)^\perp \cap R$ is a proper closed root subsystem of $R$ of maximal order.
Moreover, we can choose $\mu$ and $\nu$ such that $\mu+\nu$ lies in the root lattice.}

\proof
Let $\varpi_i$ denote the $i$th fundamental weight in the Bourbaki ordering \Bourbaki~Planches.
By definition, $\alpha_j^\vee(\varpi_i)=\delta_{ij}$, so the Dynkin diagram of the root system $\varpi_i\perp\cap R$ is obtained from that of $R$ by deleting the $i$th vertex.

All claims made below about Weyl orbits of weights $\varpi_i$ for exceptional groups were checked with \simpLie.  Note that \simpLie\ uses the Dynkin root ordering, which differs from the Bourbaki ordering for $E_n$.

For $R = A_n$, we may take 
$$\mu = -\nu = \varpi_2+\varpi_3+\cdots + \varpi_n.$$
As $\nu+\rho = \varpi_1$ is the fundamental weight associated to vertex $1$, $(\nu+\rho)^\perp\cap R$ is of type $A_{n-1}$.

For $R = B_n$, we may take
$$\mu = -\nu = \varpi_2+\varpi_3+\cdots+\varpi_n.$$
As $\nu+\rho = \varpi_1$, $(\nu+\rho)^\perp\cap R$ is of type $B_{n-1}$.

For $R = C_n$, we may take
$$\mu = -\nu = \varpi_1+\varpi_2+\cdots+\varpi_{n-2}+2\varpi_{n-1}.$$
As $\nu+\rho = \varpi_n-\varpi_{n-1}$ lies in the Weyl orbit of $\varpi_1$, $(\nu+\rho)^\perp\cap R$ is of type $C_{n-1}$.

For $R = D_n$, we may take
$$\mu = -\nu = \varpi_1+\varpi_2+\cdots+\varpi_{n-3}+2\varpi_{n-2}.$$
As $\nu+\rho = \varpi_n+\varpi_{n-1}-\varpi_{n-2}$ lies in the Weyl orbit of $\varpi_1$, $(\nu+\rho)^\perp\cap R$ is of type $D_{n-1}$.

For $R = E_6$, we may take 
$$\mu = \varpi_2 + \varpi_3+ \varpi_5 + \varpi_6,\ \nu = -\varpi_2 - 2\varpi_3 - 2\varpi_5 - \varpi_6.$$
As $\nu+\rho = \varpi_1-\varpi_3+\varpi_4-\varpi_5$ lies in the Weyl orbit of $\varpi_1$, 
$(\nu+\rho)^\perp\cap R$ is of type $D_5$.

For $R = E_7$, we may take
$$\mu = \varpi_1 + \varpi_3 + \varpi_4 + \varpi_6,\ \nu = -\varpi_1 - \varpi_3 - 2\varpi_4 - 2\varpi_6.$$
As $\nu+\rho = \varpi_2 - \varpi_4 + \varpi_5 - \varpi_6+\varpi_7$ lies in the Weyl orbit of $\varpi_7$, $(\nu+\rho)^\perp\cap R$ is of type $E_6$.

For $R = E_8$, we may take
$$\mu = \varpi_1+\varpi_3+\varpi_4+\varpi_6+\varpi_8,\ \nu = -\varpi_1-\varpi_3-2\varpi_4-2\varpi_6-2\varpi_8.$$
As $\nu+\rho = \varpi_2-\varpi_4+\varpi_5-\varpi_6+\varpi_7-\varpi_8$ lies in the Weyl orbit of $\varpi_8$, $(\nu+\rho)^\perp\cap R$ is of type $E_7$.

For $R = F_4$, we may take
$$\mu = \varpi_1+\varpi_2,\ \nu = -\varpi_1-2\varpi_2.$$
As $\nu+\rho = -\varpi _2+ \varpi_3 + \varpi_4$ lies in the Weyl orbit of $\varpi_4$, $(\nu+\rho)^\perp\cap R$ is of type $B_3$.

For $R = G_2$, we may take
$$\mu = \varpi_2,\ \nu = -2\varpi_2.$$
As $\nu+\rho = \varpi_1 - \varpi_2$ lies in the Weyl orbit of $\varpi_2$, $(\nu+\rho)^\perp\cap R$ is of type $A_1$.

Except in the exceptional cases, $\mu+\nu = 0$, so the root lattice condition is trivial.  It is likewise trivial for $E_8$, $F_4$, and $G_2$.  This leaves $E_6$ and $E_7$, where the claim can be checked by \simpLie.
\qed

We recall that the \emph{height} of a rational number $m/n$ is $\max(|m|,|n|)$, where $\gcd(m,n)=1$.

\nproc\Height{Lemma}{If $R_1$ and $R_2$ are irreducible root systems and $\mu_i$, $\nu_i$ are roots of $R_i$ for $i=1,\,2$, either $P^{R_1}_{\mu_1,\nu_1}(x)$ is a constant multiple of 
$P^{R_2}_{\mu_2,\nu_2}(x)$ or the height
$${\rm ht}\Bigl({P^{R_1}_{\mu_1,\nu_1}(p)\over P^{R_2}_{\mu_2,\nu_2}(p)}\Bigr)$$
goes to $\infty$ as $p$ increases without bound.}

\proof
More generally, if $f(x)$ is any non-constant rational function over $\Q$,
the height of $f(n)$ goes to $\infty$ as $n\to\infty$.  Indeed, writing $f(x)$ as
$a Q(x)\over  b R(x)$, where $a$ and $b$ are relatively prime integers and $Q(x)$ and $R(x)$ are monic polynomials with integer coefficients, the greatest common divisor of $Q(n)$ and $R(n)$ divides the (non-zero) resolvent of the two polynomials.  Therefore, if either $Q$ or $R$ is different from $1$, the height of $Q(n)/R(n)$ (and therefore the height of $f(n)$) goes to infinity.\qed

\sec{Efficiency of root systems}

In this section we give a definition of efficiency purely in terms of root systems and compute the efficiency of each simple root system.  It is not obvious that this coincides with the characterization \limitDef\ given in the introduction, though the equivalence can be deduced easily from Proposition 4.3 below.  

\proc{Definition}{The \emph{efficiency} $\eff(R)$ of an irreducible root system $R$ is the supremum of
$|R'|\over |R|-|R''|$, where $R'$ and $R''$ are  closed root subsystems of $R$ such that
$R'\cap R'' = \emptyset$.  
%Note that as $R'$ and $R''$ are closed root systems, disjointness is equivalent to the apparently %stronger condition
%
%$$R'\cap \Span R''=\Span R'\cap R'' = \emptyset.\eql\disjointness$$
%
The \emph{level} $\lev(R)$ is the smallest value of $|{R'}^+|$
for which $|R'|\over |R|-|R''|$ achieves $\eff(R)$.}

This definition is connected to the study of Weyl polynomials by the following proposition:

\nproc\Relevance
{Proposition}{Let $R$ be an irreducible root system with half-sum of roots $\rho$, and consider all pairs of weights $\mu$ and $\nu$ of $R$ satisfying \dominantI, \dominantII, and $\nu+\rho\neq 0$.
The maximum of 
$$\ord_{x=0}P^R_{\mu,\nu}(x)\over \deg P^R_{\mu,\nu}(x)$$
is given by $\eff(R^\vee)$, and this maximum is achieved only when $\ord_{x=0}P^R_{\mu,\nu}(x)\ge\lev(R^\vee)$.}

Note that $R = R^\vee$ except for the case $B_n/C_n$ when $n\ge 3$, and even here, as we will see, $\eff(R) = \eff(R^\vee)$ and $\lev(R) = \lev(R^\vee)$, so in fact \Relevance\ remains true when all the occurrences of $R^\vee$ are replaced with $R$.
\proof 
Let $(R^\vee)' = (\nu+\rho)^\perp \cap R^\vee$, and $(R^\vee)'' = \mu^\perp\cap R^\vee$.   As $\nu+\rho\neq 0$, $(R^\vee)''$ is a proper subsystem of $R^\vee$.  If $\mu = 0$, $P^R_{\mu,\nu}(x)$ is a non-zero constant, so the maximum will not be achieved by such a choice.  Therefore, $(R^\vee)'$ and $(R^\vee)''$ are proper closed root subsystems of $R^\vee$, disjoint by \dominantII.  A factor $\alpha^\vee(\mu)x + \alpha^\vee(\nu+\rho)$ of $P^R_{\mu,\nu}(x)$
has degree $0$ if $\alpha\in (R^\vee)''$ and $1$ otherwise.  The order of this polynomial at $0$ is $1$ if $\alpha\in (R^\vee)'$ and $0$ otherwise.   Taking the product over $R^+$, we see that
$${\ord_{x=0}P^R_{\mu,\nu}(x)\over \deg P^R_{\mu,\nu}(x)}\le\eff(R^\vee)$$
for all $\mu$ and $\nu$ satisfying the hypotheses of the proposition.  All the remains to be shown is that a pair $((R^\vee)',(R^\vee)'')$ achieving the bound $\eff(R)$ can always be found by taking a suitable choice of $\mu$ and $\nu$.  We claim that the choices in \Explicit\ do the trick.  This follows from the explicit computation of efficiencies of all irreducible root systems given by the following proposition.

\nproc\SimpleEfficiencies{Proposition}
{The efficiencies of root systems of type $A_n$, $B_n$, $C_n$, $D_n$, $E_6$, $E_7$, $E_8$, $F_4$, and $G_2$ are $n\over n+2$, $n-1\over n+1$, $n-1\over n+1$, $n-1\over n+1$, 
$10\over 17$, $3\over 5$, $63\over 117$, $3\over 7$, and $1\over 5$ respectively.
In each case the efficiency is achieved only when $R'$ is a proper closed root subsystem of maximal order in $R$.
Explicitly, the levels of these root systems are $n^2-n\over 2$, $(n-1)^2$, $(n-1)^2$, $(n-1)(n-2)$, $20$, $36$, $63$, $9$, and $1$ respectively.  
Two inequivalent irreducible root systems with the same efficiency and the same level must be of the form $B_n$ and $C_n$ for the same $n\ge 3$.}

\proof
We prove that in each case the efficiency is achieved when $R'$ is the largest proper root subsystem and $R''$ is as large as possible given $R'\cap R'' = \emptyset$.  It suffices to determine $R'_0$ and $R''_0$ satisfying these conditions and show that if $R'$ is not of maximal order, then
${|R'|\over |R|-|R''|}<{|R'_0|\over |R|-|R''_0|}$.  In particular, it suffices to prove $|R'|+|R''|\le |R'_0|+|R''_0|$.

For every closed root subsystem $S\subset R$, there exist choices of basis for $R$ and $S$ such that every positive root of $S$ is positive in $R$ and every simple root of $S$ is simple in $R$
(\Bourbaki~VI,~\S1,~Prop.~24).  Thus we can regard root subsystems as arising from subgraphs of the Dynkin diagram of $R$.  This is convenient, among other things, for computing the largest subsystems of the different root systems (which are always obtained by removing a single vertex from the Dynkin diagram).   Explicitly, the proper closed root systems of maximal order for $A_n$, $B_n$, $C_n$, $D_n$, $E_6$, $E_7$, $E_8$, $F_4$, and $G_2$ are $A_{n-1}$, $B_{n-1}$, $C_{n-1}$, $D_{n-1}$, $D_5$, $E_6$, $E_7$, $B_3/C_3$, and $A_1$ respectively.  These give the values for $\lev(R)$ stated in the proposition.

It is also worth noting that if rank $R'$ is $n-1$, and $f$ is a non-zero linear form vanishing on the hyperplane $\Span R'$, then the Coxeter number $h$ of every simple factor of $R''$ disjoint from $R'$ satisfies 
$$h\le 1+\sup_{s_1,s_2\in R\setminus R'}{f(s_1)\over f(s_2)}.\eql\levels$$
Indeed, a choice of half-plane with boundary $\Span R'$ determines a basis on $R''$.  
If $\alpha_1,\ldots,\alpha_t$ are the simple roots and $\alpha = \sum n_i\alpha_i$ the highest root of a simple factor of $R''$ with Coxeter number $h$, then 
$$\sup_{s_1,s_2\in R\setminus R'}{f(s_1)\over f(s_2)} \ge \sup_j \sum_{i=1}^t n_i {f(\alpha_i)\over f(\alpha_j)} \ge \sum _{i=1}^t n_i = h-1.$$
The value on the right hand side of \levels\ can be computed easily by restricting the adjoint representation of a group of type R to $U(1)$ times a group of type $R'$ and looking at the $U(1)$-components of the resulting direct sum.

With these general facts in hand, we now examine the cases.  

For $A_n$, we use $\varepsilon$-coordinates, in which the roots  are the vectors of length $\sqrt 2$ in $\ker\Sigma\colon\Z^{n+1}\to\Z$, $\Sigma$ denoting the sum of coordinates.  
Let $R'_0 = \varepsilon_1^\perp\cap A_n = A_{n-1}$ and 
$R''_0 = \{\pm(\varepsilon_1-\varepsilon_2)\}$.  We define an equivalence relation on 
$\{1,\ldots,n+1\}$ by setting $i\sim j$ if and only if $\varepsilon_i-\varepsilon_j\in R'\cup\{0\}$.  The equivalence classes $I_1,\ldots,I_c$ have cardinality $a_1,\ldots,a_c$ where $\sum a_k = n+1$ and $R' = \prod_{a_k > 1} A_{a_k-1}$.  Given $k\neq l$, the set $\{\varepsilon_i-\varepsilon_j\in R''\mid i\in I_k, j\in I_l\}$ has at most $\min(a_i,a_j)\le (a_i-1)(a_j-1)+1$ elements since $R'$ and $R''$ are disjoint.  Thus
$$\eqalign{|R'|+|R''| &\le \sum_k a_k^2 - a_k + \sum_{k\neq l}(a_k a_l - a_k - a_l + 2)\cr
& = (n+1)^2 - (n+1) - 2(c-1)(n+1) + 2c(c-1)\cr
& = (n+2-c)^2 + c^2 - (n+2).\cr}$$
Thus $|R'|+|R''| <n^2 - n + 2$ unless $c\le 2$ or $c\ge n$.  For $c\le 2$, $|R''|\le 2$, and for $c\ge n$, $|R'| \le 2$, so either way $|R'|+|R''|$ is bounded above by $|R'_0| + |R''_0|$.

For $B_n$ (resp.  $C_n$) we use the $\varepsilon$-coordinate system, in which the short  (resp. long) vectors are $\varepsilon_i$ (resp.  $2\varepsilon_i$).   
Let $R'_0 = \varepsilon_1^\perp\cap B_n = B_{n-1}$ (resp. $\varepsilon_1^\perp\cap C_n = C_{n-1}$)
and $R''_0 = \{\pm\varepsilon_1\}$ (resp. $\{\pm 2\varepsilon_1\}$).
Let $I'$ be the set of $i\in\{1,\ldots,n\}$ such that
$\varepsilon_i$ (resp. $2\varepsilon_i$) lies in $R'$ and $I''$ the corresponding set for $R''$.  As $R'$ and $R''$ are disjoint, $I'\cap I'' = \emptyset$, and we let $I$ denote the complement of $I'\cup I''$.
By disjointness, $\pm \varepsilon_i\pm \varepsilon_j\not\in R''$ for $i,j\in I'$ and $\pm \varepsilon_i\pm \varepsilon_j\not\in R'$ for $i,j\in I''$; the set $\{\pm \varepsilon_i\pm \varepsilon_j\mid i\in I, j\in I'\}$, which is disjoint from $R'$, meets $R''$ in at most
$2\min(|I|, |I'|)$ elements, and likewise $\{\pm \varepsilon_i\pm \varepsilon_j\mid i\in I, j\in I''\}$, which is disjoint from $R''$, meets $R'$ in at most $2\min(|I|,|I''|)$ elements.  
Finally, for each $i,j\in I$, $i\neq j$, the set $\{\pm \varepsilon_i\pm \varepsilon_j\}$ meets $R'$ in at most 2 vectors and $R''$ in at most 2 vectors.  Combining these facts, we get
$$\eqalign{|R'|+|R''|&\le (2|I'|^2 + |I|^2 - |I| + 2\min(|I|,|I''|))\cr
& \qquad\qquad\qquad\qquad +(2|I''|^2 + |I|^2 - |I| + 2\min(|I|, |I'|)) \cr
& \le (2|I'|^2 + |I|^2 - |I| + |I|+|I''|)+(2|I''|^2 + |I|^2 - |I| + |I| +  |I'|) \cr
& = 2|I|^2 + 2(|I'|+{1\over 4})^2 + 2(|I''|+{1\over 4})^2-{1\over 4}.\cr}\eql\sumR$$
As $|I| + |I'| + |I''| = n$,
for $n\ge 4$, this quantity is $\le |R'_0|+|R''_0| = 2n^2-4n+4$
unless one of $|I|$, $|I'|$, and $|I''|$ is at least $n-1$.  If $|I| = n-1$, then the first inequality of \sumR\ reads
$$|R'|+|R''|\le 2n^2-6n+8\le 2n^2-4n+4=|R'_0|+|R''_0|.$$
If $|I|=n$, then $|R'|\le n^2-n$ and $|R''|\le n^2-n$, with equality only if the systems are of type $A_{n-1}$.  If $R'$ is of type $A_{n-1}$, the right hand side of \levels\ is $3$ (resp. $2$), so assuming $n\ge 4$, $R''$ cannot also be of type $A_{n-1}$.  
Thus,
$${|R'|\over |R|-|R''|} < {n(n-1)\over 2n^2 - n(n-1)} = {n-1\over n+1} = {|R'_0|\over |R|-|R''_0|}.$$
Thus $|R'|+|R''|$ is indeed maximized when $R'$ is $B_{n-1}$ (resp. $C_{n-1}$).
The same is true for $n=2$ and $n=3$ by a simple case analysis.

For $D_n$, we again use the $\varepsilon$-coordinate system, in which the roots are the length $\sqrt 2$ vectors in $\Z^n$.   Let $R'_0 = \varepsilon_1^\perp\cap D_n = D_{n-1}$ and 
$R''_0 = \{\pm\varepsilon_1\pm\varepsilon_2\} = D_2$.
The subsystem $R'$ defines an equivalence relation $\sim$ on $\{1,2,\ldots, n\}$, where $i\sim j$ if $i=j$ or $\{\pm \varepsilon_i\pm \varepsilon_j\}\subset R'$.  As $R'$ is a closed root subsystem, there can be at most one $\sim$ equivalence class with 2 or more elements; indeed,
$$\{\pm\varepsilon_i\pm\varepsilon_k\}\subset
\Span\{\pm\varepsilon_i\pm\varepsilon_j,\pm\varepsilon_k\pm\varepsilon_l\}.$$
Let $I'$ denote this class (or the empty set) and $I''$ the corresponding set for $R''$.  As $R'\cap R'' = \emptyset$, $|I'\cap I''|\le 1$.
Let $I$ denote the complement of $I'\cup I''$.  
Reasoning as in the case $B_n/C_n$, $\pm \varepsilon_i\pm \varepsilon_j\not\in R''$ for $i,j\in I'$ and $\pm \varepsilon_i\pm \varepsilon_j\not\in R'$ for $i,j\in I''$; the set $\{\pm \varepsilon_i\pm \varepsilon_j\mid i\in I, j\in I'\}$, which is disjoint from $R'$, meets $R''$ in at most $2\min(|I|, |I'|)$ elements, and likewise $\{\pm \varepsilon_i\pm \varepsilon_j\mid i\in I, j\in I''\}$, which is disjoint from $R''$, meets $R'$ in at most $2\min(|I|,|I''|)$ elements.  Again  for each $i,j\in I$, $i\neq j$, the set $\{\pm \varepsilon_i\pm \varepsilon_j\}$ meets $R'$ in at most 2 vectors and $R''$ in at most 2 vectors, so
$$\eqalign{|R'|+|R''|&\le (2|I'|^2- 2|I'| + |I|^2 - |I| + 2\min(|I|,|I''|))\cr
& \qquad\qquad\qquad\qquad +(2|I''|^2 - 2 |I''| + |I|^2 - |I| + 2\min(|I|, |I'|)) \cr
& \le  2|I|^2 + 2|I'|^2 - |I'| + 2|I''|^2 - |I''|\cr
& = 2|I|^2 + 2(|I'|-{1\over 4})^2 + 2(|I''|-{1\over 4})^2-{1\over 4}.\cr}$$
Moreover, $|I'|,|I''|\le n-1$, $|I|+|I'|\le n$, $|I|+|I''|\le n$, and $|I|+|I'|+|I''|\le n+1$.
If
$$|R'| + |R''| > |R'_0|+|R''_0| = 2n^2-6n+8,$$
either $|I|$, $|I'|$, or $|I''|$ is at least $n-1$.
Now $|I|\ge n-1$ implies $I' = I'' = \emptyset$, so either $R' = A_{n-1}$ or $|R'|\le (n-1)(n-2)$, and likewise for $R''$.  If $|R'|\le (n-1)(n-2)$ and $|R''|\le (n-1)(n-2)$, then ${|R'|\over |R| - |R''|}\le {n-2\over n+1} <{ |R'_0|\over |R|-|R''_0|}$.  If $R' = A_{n-1}$, the upper bound \levels\ is 2, so 
$R''$ is of the form $A_1^r$ (and vice versa).  If $R'$ or $R''$ has $r$ mutually orthogonal roots,
since they are closed root subsystems,
the Dynkin diagram of $D_n$ must have $r$ vertices, no two of which share an edge, so $r\le n/2+1$, and
$$|R'|+|R''|\le n(n-1)+n+2 \le 2n^2 - 6n + 8 = |R'_0|+|R''_0|,$$
as long as $n\ge 5$.  For $n=4$, we check by hand that there do not exist 3 mutually orthogonal roots in $D_4\setminus A_3$, so $(D_3,D_2)$ is still optimal.

For $E_8$, $|R| = 240$.  The largest closed subsystem, $E_7$, has 126 roots, and it can be realized as the set of roots orthogonal to a given root $\alpha$;  $D_7$ is next with $84$.   
As ${126\over 240-6}  >  {84\over 240-84}$, it suffices to prove that the largest closed root subsystem of $E_8$ disjoint from $E_7 = \alpha^\perp\cap E_8$ has order 6.
The orthogonal projection of any root onto $\R\alpha$ is in $\{-\alpha,-\hf\alpha,0,\hf\alpha,\alpha\}$, so the right hand side of \levels\ is $3$.  Thus all factors of $R''$ are of type $A_1$ or $A_2$.  Moreover, a factor of type $A_2$ must contain $\alpha$ itself (which is the only root whose projection onto $\R\alpha$ equals $\alpha$), so if $R''$ contains any factor of the form $A_2$, then $R'' = A_2$.
Otherwise, $R'' = A_1^n$.   We look for the largest value such that $A_1^n$ can be embedded as a closed subsystem in $E_8$ with no root orthogonal to $\alpha$.  Equivalently, 
we fix an embedding of $A_1^n$ in $E_8$ and ask whether any root of $E_8$ fails to be orthogonal to every root of $A_1^n$.  Equivalently, we ask whether the restriction of the adjoint representation of $E_8$ to $A_1^n$ has an irreducible summand which is non-trivial on each factor $A_1$.

To construct mutually orthogonal roots in $E_8$, we start with a root $\beta_1$, choose $\beta_2\in \beta_1^\perp\cap E_8 = E_7$, choose $\beta_3\in \beta_2^\perp\cap E_7 = D_6$, choose $\beta_4\in \beta_3^\perp\cap D_6 = A_1\times D_4$, and so on.  All choices are equivalent until $\beta_4$, which can be chosen in $A_1$ or in $D_4$.  In the first case, $A_1^4\subset D_4^\perp\cap E_8 = D_4$, so the system is not closed.  We use branching rules \simpLie, starting with the adjoint representation of $E_8$, restricting to $A_1\times E_7$, taking only representations of $E_7$ which appear opposite a non-trivial factor in $A_1$, and so on.  We get the (56-dimensional) representation $V_{\varpi_7}$ of $E_7$, and then the (12-dimensional) representation $V_{\varpi_1}$ of $D_6$, then a representation of $A_1\times D_4$ which is 2-dimensional and therefore trivial on the $D_4$ factor.  We conclude that $A_1^3$ can be realized (and it is indeed a closed subsystem of $E_8$ since no simply laced rank 3 system contains $A_1^3$ as a proper subsystem), but $A_1^4$ cannot.

For $E_7$, $|R|=126$.   The largest subsystem, $E_6$, has 72 roots; $D_6$ comes next, with 60, and then $A_6$, with 42.  As ${72\over 126-6} > {42\over 126-42}$, we need only consider
the cases $R' = E_6$ and $R' = D_6$.  Writing $E_6$ as the orthogonal complement of the fundamental weight $\varpi_7$ the right hand side of \levels\ is $2$, so we need only consider 
$R''  = A_1^r$.   We proceed as before, with $\beta_1\in E_7$, $\beta_2\in \beta_1^\perp\cap E_7 = D_6$, and $\beta_3\in \beta_2^\perp\cap D_6 = A_1\times D_4$.  Applying the branching rules for $V_{\varpi_7}$, which we have already described, we see that if we choose $\beta_3$ to be in the factor $A_1$, the restriction of $V_{\varpi_7}$ to $A_1^3$ has a summand non-trivial on every factor, so we can take $R''=A_1^3$.   Any fourth root $\beta_4$ mutually orthogonal to the first three must lie in the $D_4$ factor, and we know from the $E_8$ analysis that $\varpi_7$ cannot be non-orthogonal to $\beta_1$, $\beta_2$, and such a $\beta_4$.
When $R'=D_6 = \alpha^\perp\cap E_7$, the argument for $E_8$ shows that either $R'' = A_2$ or $R'' = A_1^n$ where $n\le 7$ for dimension reasons.  As ${72\over 126-6} > {60\over 126-14}$, we see that $(E_6,A_1^3)$ is optimal.

For $E_6$, $|R| = 72$.  The largest subsystem, $D_5$, has 40 roots, the next largest, $A_5$, has 30, and the third largest, $D_4$, has 24.   As ${40\over 72-4} > {24\over 72-24}$, we need only consider
$D_5 = \varpi_1^\perp\cap E_6$ and $A_5 = \alpha^\perp\cap E_6$, where $\alpha$ is a root.
For $R' = D_5$, the right hand side of \levels\ is $2$, so $R'' = A_1^n$.   Choose $\beta_1\in E_6$, $\beta_2\in \beta_1^\perp\cap E_6 = A_5$, and $\beta_3\in\beta_2^\perp\cap A_5 = A_3$. 
In the restriction of $V_{\varpi_1}$ to $A_1\times A_5$, the only summand non-trivial on $A_1$ is $6$-dimensional on $A_5$; its restriction to $A_1\times A_3$ is the a direct sum of a representation trivial on $A_1$ and a representation trivial on $A_3$.  It follows that the restriction of $V_{\varpi_1}$ to $A_1^2$ has a summand non-trivial on both factors but the restriction to $A_1^3$ has no such summand.
For $R' = A_5$, the right hand side of \levels\ is $3$, so as in the $E_8$ case, $R'' = A_2$ or $R'' = A_1^n$, and of course $n\le 5$.  As ${40\over 72-4} > {30\over 72-10}$, $(R',R'') = (D_5,A_1^2)$ is optimal.

For $F_4$, $|R| = 48$, and the largest proper closed root subsystems are obtained by deleting an endpoint vertex from the Dynkin diagram.  They are of type $B_3$ and $C_3$, and are both of order $9$; in each case, $\Span R' = \alpha^\perp\cap R$ for a root $\alpha$, short or long respectively.  Either way, the bound on the Coxeter number of components of $R''$ is $3$, so we need only consider $A_1$, $A_2$, $A_1\times A_1$, and $A_2\times A_1$ as possibilities for $R''$.
In each case, $A_2$ can be realized; for example, if $\alpha = \varpi_4$ with respect to the standard coordinates on $F_4$, $R''$ can be taken to have basis $\{2\varpi_4-\varpi_3,\varpi_3-\varpi_4\}$.  In each case, $\alpha$ or $-\alpha$ is the highest vector of the $A_2$, so $A_2\times A_1$ cannot be disjoint from $\Span R'$ (since $\pm \alpha$ belong to both).  Finally, if $R'$ is not of type $B_3$ or $C_3$, it has at most 8 roots, and $8+8 < 18+6$.

The case $G_2$ is trivial since any proper root subsystem has rank 1 and therefore order 2.

\qed

\proc{Definition}{We say an irreducible root system $S$ is \emph{greater than or equal to}
an irreducible root system $T$ (and write $S\ge T$) if and only if $\eff(S) > \eff(T)$ or $\eff(S) = \eff(T)$ and $\lev(S)\le \lev(T)$.  If $S\ge T$ but $T\ngeq S$, we say $S$ \emph{is strictly greater than} $T$.}

An immediate consequence of \SimpleEfficiencies\ is that $S\ge T$ and $T\ge S$ implies $S = T$ or $S = T^\vee$.  

\sec{The main theorem}

\nproc\Quadratic
{Lemma}{Let $R$ be an irreducible root system spanning a Euclidean space $E$ and $Q$ a quadratic form on $E$ such that $Q(\alpha) = 0$ for all $\alpha\in R$.  Then $Q = 0$.}
\proof
Let $B(x,y) = Q(x+y)-Q(x)-Q(y)$ denote the associated bilinear form.  If $\alpha$ and $\beta$ are roots of $R$ making an obtuse angle, then $\alpha+\beta$ is again a root, so $B(\alpha,\beta) = 0$.  Likewise, if $\alpha$ and $\beta$ make an acute angle, $\alpha-\beta$ is a root, so again $B(\alpha,\beta)=0$. The roots in $R\setminus \alpha^\perp$ span $E$ by \Spanning.  Therefore, $B(\alpha,E) = 0$.  As every $\alpha\in R$ lies in the null space of $B$, is follows that $B$ and therefore $Q$ must be $0$.\qed

\nproc\Ratio
{Proposition}
{If $R_1$ and $R_2$ are irreducible root systems with $R_1\ge R_2$, then there exist weights $\mu_1$ and $\nu_1$ and an integer $N$ such that:
\medskip
\item{1)}$\mu_1$ and $\nu_1$ satisfy \dominantI\ and \dominantII;
\item{2)}$\nu_1+\rho_1\neq 0$;
\item{3)}$\mu_1+\nu_1$ lies in the root lattice of $R_1$;
\item{4)}${\ord_{x=0}P^{R_1}_{\mu_1,\nu_1}(x)\over \deg P^{R_1}_{\mu_1,\nu_1}(x)} = \eff(R_1)$;
\item{5)}If $\mu_2$ and $\nu_2$ are weights of $R_2$ and $c$ is a positive constant such that 
$$P^{R_2}_{\mu_2,\nu_2}(x)=c P^{R_1}_{\mu_1,\nu_1}(x),\eql\scalar$$
then 
\itemitem{a)}$\eff(R_1) = \eff(R_2)$ and $\lev(R_1) = \lev(R_2)$.
\itemitem{b)}$c \ge 1$, with equality only if $R_1 = R_2$;
\itemitem{c)}If $p > N$ and $\ord_p c > 0$, then $p$ divides $\nu_2+\rho_2$, i.e., $\nu_2$ is not alowable at $p$.}

\proof 

First we observe that condition (4) implies that $P^{R_1}_{\mu_1,\nu_1}(x)$ is not a monomial in $x$, so that \scalar\ implies the same for $P^{R_2}_{\mu_2,\nu_2}(x)$.  This implies that $\nu_2+\rho_2\neq 0$, so 
$${\ord_{x=0} P^{R_2}_{\mu_2,\nu_2}(x)\over \deg P^{R_2}_{\mu_2,\nu_2}(x)} \le \eff(R_2).$$
The hypothesis $\eff(R_1)\ge\eff(R_2)$ then implies $\eff(R_1) = \eff(R_2)$.  \Relevance\ and equation \scalar\ then imply 
$$\lev(R_2)\le \ord_{x=0}P^{R_2}_{\mu_2,\nu_2}(x) =  \ord_{x=0}P^{R_1}_{\mu_1,\nu_1}(x) = \lev(R_1),$$
so $\lev(R_2) = \lev(R_1)$.

For (b), we choose $\mu_1$ and $\nu_1$ as in \Explicit.   We note that for any integer $n$, whether or not $\lambda_2 = n\mu_2+\nu$ is dominant, $P^{R_2}_{\mu_2,\nu_2}(n)$ is either zero or $\pm1$ times the dimension of an irreducible representation of a Lie group of type $R_2$.  Indeed, suppose $\alpha^\vee(\lambda_2+\rho_2) \neq 0$ for all $\alpha$.  Then defining $w$ in the Weyl group of $R_2$ as the unique element such that $w(\lambda_2+\rho_2)$ is dominant, and setting $\lambda'_2 = w(\lambda_2+\rho_2)-\rho_2$, the weight $\lambda'_2$ is dominant.  Thus,
$$\eqalign{\dim V_{\lambda'_2}^2 
&= (-1)^{|R_2^+|}{\prod_{\alpha\in R}\alpha^\vee(\lambda'_2+\rho_2)\over \prod_{\alpha\in R}\alpha^\vee(\rho_2)}\cr
&=(-1)^{|R_2^+|}{\prod_{\alpha\in R}\alpha^\vee(\lambda_2+\rho_2)\over \prod_{\alpha\in R}\alpha^\vee(\rho_2)}
=  P^{R_2}_{\mu_2,\nu_2}(n)^2.\cr}$$
For $R_1$ of type $A_n$, $B_n$, $C_n$, or $D_n$, we have $P^{R_1}_{\mu_1,\nu_1}(1) = 1$,
so  $P^{R_2}_{\mu_2,\nu_2}(1) = c$ implies $c\in\Z$ and therefore $c\ge 1$.   The same is true for $G_2$, since $P^{R_1}_{\mu_1,\nu_1}(2) = 1$.  For $F_4$, 
$$\gcd(P^{R_1}_{\mu_1,\nu_1}(2),P^{R_1}_{\mu_1,\nu_1}(3)) = \gcd(52,340119) = 1,$$
so again $c$ must be an integer.
For $E_6$, we have 
$$\gcd(P^{R_1}_{\mu_1,\nu_1}(2),\!P^{R_1}_{\mu_1,\nu_1}(3),\!P^{R_1}_{\mu_1,\nu_1}(4)) \!=\!
\gcd(1728, 3171108447, 71292900343808) \!=\!1,$$
so  again $c\in\Z$.  For $E_7$, we have
$$\gcd(P^{R_1}_{\mu_1,\nu_1}(2),P^{R_1}_{\mu_1,\nu_1}(3)) = \gcd(573440,33940969546604175) = 5,$$
so $5c\in\Z$.   Thus, $|P^{R_2}_{\mu_2,\nu_2}(2)|$ is a positive integer multiple of ${573440\over 5} = 114688$.   However, the lowest degree of any irreducible representation of a Lie group of type $E_7$ which is divisible by $114688$ is $573440$ itself.  For $E_8$,
$$\eqalign{\gcd(P^{R_1}_{\mu_1,\nu_1}(2),P^{R_1}_{\mu_1,\nu_1}(3)) &= \gcd(4096000,2665014302693985712862760000) \cr &= 8000,\cr}$$
and the smallest dimension of  any irreducible representation of a Lie group of type $E_8$ which is divisible by $512 = {4096000\over 8000}$ is $4096000$.  

To finish proving (b), it suffices to show that if $R_1$ is $B_n$ or $C_n$, $n\ge 3$, $R_2 = R_1^\vee$,
$\mu_1$ and $\nu_1$ are defined as in \Explicit, and $\mu_2$ and $\nu_2$ are any weights of $R_2$, then $P^{R_1}_{\mu_1,\nu_1}\neq \pm P^{R_2}_{\mu_2,\nu_2}$.   If $\mu$ and $\nu$ are weights of $B_n$ such that $\ord_{x=0}P^{B_n}_{\mu,\nu}(x)$ achieves the maximum $(n-1)^2$, then $\nu+\rho$ must be a scalar multiple $d\varepsilon_i$ of an $\varepsilon$-basis vector, and writing
$\mu = (a_1,\ldots,a_n)$, either $a_i = \pm a_j$ for some $j\neq i$ or some $a_j = 0$.  We can always replace $\nu$ by any linked weight $w(\nu+\rho)-\rho$ with the effect of changing $P^{B_n}_{\mu,\nu}$ by a factor of $\pm 1$, so without loss of generality we may assume $\mu$ is dominant, which means we may assume either $a_1 > a_2 > \cdots >  a_{n-1}  > a_n = 0$ (and $i = n$) or, renumbering,  
$a_1 >  \cdots > a_i = a_i > \cdots > a_{n-1} > 0$.   In the first case,
$$x^{-(n-1)^2}P^{B_n}_{\mu,\nu}(x)
= {d 2^n F(a_1,\ldots,a_{n-1}) 
	\prod_{i=1}^{n-1}(a_i x - d)(a_i x + d)\over \Delta_n}\eql\CaseI
$$
where
$$\Delta_n = \prod_{\alpha\in B_n}\alpha^\vee(\rho_{B_n}) 
= \prod_{\alpha\in C_n}\alpha^\vee(\rho_{C_n})
%\prod_I{i=1}^n (2i-1) \prod_{i=1}^{n-1}\prod{j=i+1}^n (i+j)(i+j-1),
$$
and
$$F(a_1,\ldots,a_{n-1}) = \prod_{i=1}^{n-1} a_i \prod_{i=1}^{n-2}\prod_{j=i+1}^{n-1}(a_i^2-a_j^2).$$
In the second case, $x^{-(n-1)^2}P^{B_n}_{\mu,\nu}(x)$ is
$${d 2^n F(a_1,\ldots,a_{n-1})(a_i x + d)
	(2a_i x+d)\prod_{j\neq i} ((a_i + a_j)x + d)((a_j - a_i)x + d)\over \Delta_n}\eql\CaseII
$$
If $A$ is the number of positive roots of  $x^{-(n-1)^2}P^{B_n}_{\mu,\nu}(x)$ and $B$ is the number of negative roots, then \CaseI\ applies if and only if $A = B$, and otherwise \CaseII\ applies
with $i = \sup(A,B) - (n-1)$.
In either case, the roots determine the vector $(a_1,\ldots,a_{n-1},d)$ up to scalar multiplication by a positive rational number and therefore the polynomial itself up to multiplication by an element of $(\Q^*)^{n^2}$.

The case analysis for $P^{C_n}_{\mu,\nu}(x)$ is exactly the same, and the formulas \CaseI\ and \CaseII\ the same except without the factor of $2^n$.  As $2^n\not\in (\Q^*)^{n^2}$, 
$P^{R_1}_{\mu_1,\nu_1}\neq P^{R_2}_{\mu_2,\nu_2}$.

Given two pairs $(\mu_1,\nu_1)$ and $(\mu_2,\nu_2)$ of weights of an irreducible root system $R$,  $i=1,\,2$, we have maps $f_i\colon R^+\to \Q^2$ defined by
$$f_i(\alpha) = (\alpha^\vee(\mu_i),\alpha^\vee(\nu_i)).$$
Suppose that $P^R_{\mu_2,\nu_2}(x)=c P^R_{\mu_1,\nu_1}(x)\neq 0$.   Then $(0,0)\not\in f_i(R^+)$ and for some permutation $\sigma\colon R^+\to R^+$, $f_2(\alpha)$ is a scalar multiple of $f_1(\alpha)$ for all $\alpha$.  As there are only finitely many possibilities for $\sigma$, by making $N$ large enough we may assume without loss of generality that $\sigma$ is the identity.  Thus,
$$\alpha^\vee(\mu_1)\alpha^\vee(\nu_2+\rho_2) - \alpha^\vee(\mu_2)\alpha^\vee(\nu_1+\rho_1) = 0$$
for all $\alpha^\vee\in R^\vee$.   Regarding the above expression as a quadratic form on $R^\vee$, by \Quadratic, it is trivial, so $(\mu_2,\nu_2+\rho_2) = C(\mu_1,\nu_1+\rho_1)$ for some non-zero rational number $C$.   Therefore $c = C^{|R^+|}$, and claim (c) follows.
\qed

\nproc\Finiteness
{Proposition}{Let $l\in\N$, $e\in\Q$, and $G$ a simply connected almost simple Lie group such that
for every almost simple factor $G_i$ of $G$ with root system $R_i$, either $\eff(R_i) < e$ or $\eff(R_i)=e$ and $\lev(R_i) \ge l$.  Let $A(G)$ denote the set of pairs $(p,\lambda)$, where $p$ is a prime and $\lambda$ is a dominant weight of $G$ allowable at $p$.  Then
$$\inf_{(p,\lambda)\in A(G)} e\log\dim V_\lambda - (\log p)\ord_p\dim V_\lambda > -\infty.\eql\lowerBound$$
Moreover, for every $C\in\R$, there exists a finite set
$$L_{G,e,l,C} = \{(\mu,\nu)\mid\mu,\nu\in\Lambda_G\}$$
such that for all $(p,\lambda)\in A(G)$ such that $\ord_p \dim V_\lambda \le l$ and
$$e\log\dim V_\lambda - (\log p)\ord_p\dim V_\lambda \le C,\eql\upperBound$$
there exists $(\mu,\nu)\in L_{G,e,l,C} $ such that $\lambda = p\mu+\nu$ and
either $\mu = 0$, or $\mu$ is a non-zero weight of an
almost simple factor $G_i$ of $G$ with $\eff(R_i) = e$ and $\lev(R_i) = l$.}

\proof 
Suppose first that $G$ is almost simple with root system $R$ and Coxeter number $h$.
We take advantage of the equalities $|R| = |R^\vee|$, $\eff(R) = \eff(R^\vee)$, and $\lev(R) = \lev(R^\vee)$ to avoid burdening our notation unnecessarily with $^\vee$ superscripts.
Define
$$S_i = \{\alpha^\vee\in R^\vee:\; p^i\vert \alpha^\vee(\lambda+\rho)\},\ R'_i = R^\vee\cap\Span S_i\eql\subsystem$$
Thus,
$$R^\vee\supset R'_1\supset R'_2\supset\cdots\supset R'_k\supsetneqq R'_{k+1}=\emptyset.$$
An examination of the list of all equal-rank root subsystems of all irreducible root systems shows that $\alpha\in R'_i$ implies either
$2\alpha$ or $3\alpha$ lies in the abelian group generated by $S_i$.
As $\lambda$ is allowable at $p$, $R'_i$ is a proper subsystem of $R^\vee$ except possibly when $i=1$ and $p\le 3$.
Let $R''_i$ denote the closed root subsystem of $R^\vee$
generated by the simple roots $\alpha^\vee$ such that
$$0 < \alpha^\vee(\lambda+\rho) < {p^i\over  3h}.\eql\compineq$$
For any root in $R''_i$, $0 < \alpha^\vee(\lambda+\rho) < {p^i\over 3}$, so neither $2\alpha^\vee$ nor $3\alpha^\vee$ lies in the abelian group generated by $S_i$, 
which means $R'_i\cap R''_i = \emptyset$.  We have
$$\ord_p \dim V_\lambda \le \sum_{i=1} ^k|(R'_i)^+|\eql\totalp$$
with equality if $p\ge 5$, and
$$\log \dim V_\lambda \ge -|R^+| \log 3h+ \log p\sum_{i=1}^k (|R^+| - |(R''_i)^+|).\eql\logdim$$
Except possibly when $i=1$ and $p\le 3$,
$$\eff(R)(|R^+| - |(R''_i)^+|) \ge |(R'_i)^+|,\eql\defineq$$
so
$$\eff(R)\log p\sum_{i=1}^k (|R^+| - |(R''_i)^+|)\ge -\eff(R)|R^+|\log 3+\log p\sum_{i=1} ^k|(R'_i)^+|.\eql\efflog$$
By \logdim, \efflog, and \totalp,
$$\eqalign{e\log\dim V_\lambda -  \log p\;\ord_p \dim V_\lambda 
&\ge  \eff(R) \log \dim V_\lambda -  \log p\;\ord_p \dim V_\lambda\cr
& \ge -\eff(R) |R^+|\log 9h.\cr}$$
This gives \lowerBound.  

If $e > \eff(R)$, then 
$$
e\log\dim V_\lambda -  \log p\;\ord_p \dim V_\lambda\ge (e-\eff(R))\log\dim V_\lambda
-\eff(R) |R^+|\log 9h
$$
implies that the left hand side is $\ge C$ for all but finitely many weights $\lambda$, independent of $p$.
We can therefore define $L_{G,e,l,C}$ to be the set of pairs $(0,\nu)$, as $\nu$ ranges over all such $\lambda$.    If $l < \lev(R)$, then for $p\ge 5$, \totalp\ is an equality, and so $|(R_i')^+| < \lev(R)$ for 
$1\le i\le k$.  It follows that \defineq\ is strict and that 
$$\displaylines{\quad
\eff(R)\log\dim V_\lambda -  \log p\;\ord_p \dim V_\lambda  \hfill\cr
\hfill\ge \bigl\{\eff(R)(|R^+| - |(R''_i)^+|) - |(R'_i)^+|\bigr\}\log p-\eff(R) |R^+|\log 9h > C\quad\cr}
$$
for $p\gg 0$.  We therefore obtain an upper bound on $p$, which then gives an upper bound on $\dim V_\lambda$ and consequently a finite set of possible $\lambda$, and we can set 
$L_{G,e,l,C}$ to be a set of pairs $(0,\lambda)$ as before.  The only way in which $\dim V_\lambda$ could fail to be bounded by the condition \upperBound\ is if $\eff(R) = e$ and $\lev(R) = l$ and moreover
$|(R'_i)^+| = \lev(R)$ for all $i\le k$.  For $p\ge 5$, equality holds in \totalp, and this last condition implies $k=1$.

To prove that $\lambda$ can be written $p\mu+\nu$, where $\mu$ and $\nu$ are drawn from finite sets independent of $p$, we may assume that we have fixed $S = R'_1$ and $T = R''_1$.  The Weyl formula gives
$$\eqalign{\log \dim V_\lambda &\ge (|S^+|)\log p + 
\sum_{\alpha\in S^+}\log{\alpha^\vee(\lambda+\rho)\over p}+\sum_{\alpha\in T^+}\log \alpha^\vee(\lambda+\rho)\cr
& \qquad\qquad\qquad
+ (|R^+| - |S^+|-|T^+|)\log {p\over 3h} - \sum_{\alpha\in R^+}\log  \alpha^\vee(\rho)\cr
& = (|R^+|-|T^+|)\log p + 
\sum_{\alpha\in S^+}\log{\alpha^\vee(\lambda+\rho)\over p}+\sum_{\alpha\in T^+}\log \alpha^\vee(\lambda+\rho)\cr
&  \qquad\qquad\qquad
- (|R^+| - |S^+|-|T^+|)\log 3h - \sum_{\alpha\in R^+}\log  \alpha^\vee,\cr}
$$
so there are only finitely many possible values for the integers $p^{-1}\alpha^\vee(\lambda+\rho)$ for $\alpha\in S^+$ and $\alpha^\vee(\lambda+\rho)$ for $\alpha\in T^+$ which will bound
$$\log \dim V_\lambda - (|R^+|-|T^+|)\log p = {\eff(R)\log\dim V_\lambda -  \log p\;\ord_p \dim V_\lambda\over\eff(R)}$$
above.  For every $\beta^\vee\in \Z R^\vee\cap \Span S\cap \Span T$, $\beta^\vee(\lambda+\rho)$ is both divisible by $p$ and bounded absolutely.  Fixing a basis of $\Span S\cap \Span T$ in the root lattice of $R^\vee$,
we conclude that for $p\gg 0$, $\beta^\vee(\lambda+\rho)=0$ for all $\beta^\vee\in\Span S\cap \Span T$.  
As $S^\vee$ and $T^\vee$ generate the root lattice, 
there exists a unique weight $\nu$ such that $\alpha^\vee(\nu+\rho) = 0$ for all $\alpha^\vee\in S$ and $\alpha^\vee(\nu+\rho) = \alpha^\vee(\lambda+\rho)$ for all $\alpha^\vee\in T$.  
Letting $\mu = (\lambda-\nu)/p$, we see that $\alpha^\vee(\mu)$ is an integer for all $\alpha^\vee\in S\cup T$ and therefore for all roots in $R^\vee$; it follows that $\mu$ is also a weight, and as $\mu$ and $\nu$ are determined independently of $p$ by the 
values $p^{-1}\alpha^\vee(\lambda+\rho)$ for $\alpha^\vee\in S^+$ and $\alpha^\vee(\lambda+\rho)$ for $\alpha^\vee\in T$, there are only finitely many possibilities for $(\mu,\nu)$.  Thus the proposition is true when $G$ is almost simple.  Note that in every case where $\mu\neq 0$ arises, $\ord_p \dim V_\lambda = l$, the maximum allowable value.

Now consider a general $G = G_1\times\cdots\times G_n$.
Let
$$B_i = \inf_{(p,\lambda)\in A(G_i)}
 e\log\dim V_\lambda - (\log p)\ord_p\dim V_\lambda$$
Then letting $\lambda_i$ denote the $i$th component of $\lambda$,
$$\eqalign{e\log\dim V_\lambda - (\log p)\ord_p\dim V_\lambda 
&= \sum_i e\log\dim V_{\lambda_i} - (\log p)\ord_p\dim V_{\lambda_i}\cr
& \ge \sum B_i,\cr}$$
so $B = \sum_i B_i$ is a lower bound for \lowerBound.
On the other hand, 
$$e\log\dim V_\lambda - (\log p)\ord_p\dim V_\lambda \le C$$
implies
$$e\log\dim V_{\lambda_i} - (\log p)\ord_p\dim V_{\lambda_i} \le C - \sum_{j\neq i} B_j, = C-B+B_i$$
so there exist sets $L_{G_i,e,l,C-B+B_i}$ such that each $\lambda_i = p\mu_i+\nu_i$ for some $(\mu_i,\nu_i)\in L_{G_i,e,l,C-B+B_i}$.  Letting
$$L_{G,e,l,C} = \Bigl\{\Bigl(\sum_i \mu_i,\sum_i \nu_i\Bigr)\Bigm|(\mu_i,\nu_i)\in L_{G_i,e,l,C-B+B_i}\Bigr\},$$
we obtain a finite set of pairs $(\mu,\nu)$ of the desired kind.  In particular, for each $\lambda$ satisfying the condition $\ord_p V_\lambda\le l$, there can be at most one coordinate $i$ with $\mu_i\neq 0$ because
$$\ord_p \dim V_{p\mu_i+\nu_i} = l$$
for each non-zero $\mu_i$ and each $p\gg 0$.

\qed

\proc{Theorem}{The zeta-function of a compact semisimple Lie group $G$ determines the
Lie algebra of the group up to isomorphism.}

\proof 
Let $G^1$ and $G^2$ be compact semisimple Lie groups with $\zeta_{G^1}(s) = \zeta_{G^2}(s)$ but with non-isomorphic Lie algebras.   By \Sizes, the number of factors of type $A_1$ in $G^1$ and $G^2$ are the same.  
Let $R^i_1,\,R^i_2,\ldots$ be the irreducible factors of the root system $R^i$ of $G^i$.   If $G^1$ and $G^2$ have different Lie algebras, by a suitable renumbering, we may assume that for some $k\ge 0$
$$R^1_1\ge R^1_2\ge\cdots,\ R^2_1\ge R^2_2\ge\cdots,\ R^1_i = R^2_i\ \forall i\le k,$$
and there is no factor in $R^2_{k+1},R^2_{k+2},\ldots$ which equals $R^1_{k+1}$ or is strictly greater than it.  In particular, $R^1$ has at least $k+1$ irreducible factors (though $R^2$ might have only $k$.)  Moreover, $R^1_{k+1}$ is not of type $A_1$.
Under these conditions, 
we will find an integer $n$ such that the $n$th coefficient of $\zeta^*_{G^1}(s)$ is positive and the $n$th coefficient of 
$\zeta^*_{G^2}(s)$ is zero, which by \Sizes\ implies $\zeta_{G^1}(s)\neq\zeta_{G^2}(s)$.
Without loss of generality we may assume that $G^1$ is adjoint and $G^2$ is simply connected, so both are products of almost simple factors: $G^i=\prod_j G^i_j$.

In fact, we will prove the following stronger statement: under these conditions,  if $r$ is a fixed positive integer, then there exist primes $p_1\gg\cdots\gg p_k$, all congruent to 1 (mod $N_{G^1} = N_{G^2}$) and an allowable weight $\lambda^1$ of $G^1$ such that 
$$\ord_{p_i} \dim V_{\lambda^1}=\lev(R^1_i),\ 1\le i\le k+1\eql\exactDivis$$
and there is no allowable weight $\lambda^2$ of $G^2$ such that the following conditions hold:
$$\ord_{p_i}\dim V_{\lambda^2} = \lev(R^1_i),\ 1\le i\le k+1;\eql\inductionI$$
$$\dim V_{\lambda^2} \mid r\dim V_{\lambda^1};\eql\inductionII$$
$$\dim V_{\lambda^2}\le \dim V_{\lambda^1}.\eql\inductionIII$$
Explicitly, we may take $\lambda^1=\lambda^1_1+\cdots+\lambda^1_{k+1}$, where $\lambda^1_i = p_i \mu^1_i + \nu^1_i$, and $\mu^1_i$ and $\nu^1_i$ are
defined by \Explicit\ as weights of the factor $\tilde G^1_i$ of $\tilde G^1$.  Note that the congruence condition on $p_i$ guarantees that each $\lambda^1_i$ lies in the root lattice of $\tilde G^1_i$ and therefore defines a representation of $G^1_i$.  Note also that we take all primes $p_i$ in the construction to be larger than the bound $N$ in \Ratio\ (which depends only on the root systems of $G^1$ and $G^2$.)  We also choose each $p_i$ large enough that for each $j\le k+1$, $p_i>P^{R^1_j}_{\mu^1_j,\nu^1_j}(2)$; again, this condition depends only on $G^1$.  By choosing each $p_j$ to be congruent to $2$ modulo a sufficiently high power of $p_{j+1}\cdots p_k p_{k+1}$, we can guarantee that $p_i\ndiv \dim V_{\lambda^1_j}$ for $i>j$, and therefore, we can guarantee \exactDivis.

Suppose the claim does not hold for $k=0$.  If $R^1_1 = A_1$, we are done.
Otherwise, $P^{R^1_1}_{\mu_1^1,\nu_1^1}(x)$ is not constant.  As 
$$
\eff(R^1_1) = {\ord_{x=0}P^{R^1_1}_{\mu_1^1,\nu_1^1}(x)\over \deg P^{R^1_1}_{\mu_1^1,\nu_1^1}(x)},$$
we have
$$\eff(R^1_1)\log\dim V_{\lambda^1_1} \le \lev(R^1_1)\log p_1 + C,\eql\notTooBig$$
for some constant $C$ depending only on $R^1$.  By \inductionII,
$$\eff(R^1_1)\log\dim V_{\lambda^2}\le  \lev(R^1_1)\log p_1 + \eff(R^1_1)\log r + C.$$
Applying \Finiteness\ with $e=\eff(R^1_1)$, $l=\lev(R^1_1)$, $G=G^2$, $p=p_1$, and 
$\lambda^2$ an allowable weight of $G^2$ satisfying \inductionI--\inductionIII, we deduce that $\lambda^2$ must be of the form $p_1\mu+\nu$, where $\mu$ and $\nu$ belong to finite sets independent of $p_1$.  As $P^{R^1_1}_{\mu_1^1,\nu_1^1}(x)$ is not constant, $\dim V_{\lambda^1_1}\to\infty$ as $p_1\to\infty$, so for $p_1\gg 0$, $\mu\neq 0$.  Again by \Finiteness, $\mu$ must be a weight of an almost simple factor $G^2_1$ of $G^2$, and $R^2_1$ must have the same efficiency and level as $R^1_1$. 
(By hypothesis, $R^2_1$ and $R^1_1$ cannot be equal, so at this point we must be in the $B_n/C_n$ case, though we do not use this fact explicitly.)
Writing $G^2 = G^2_1\times H^2$, we can decompose $\lambda^2 = \lambda^2_1+(\lambda^2-\lambda^2_1)$, where the second summand is a weight of $H^2$, and $\lambda^2_1 = p_1\mu^2_1+\nu^2_1$, where $\mu^2_1$ and $\nu^2_1$ belong to finite sets independent of $p_1$.
By \lowerBound\ and \notTooBig,
the fraction $\dim V_{\lambda^1_1}\over \dim V_{\lambda^2_1}$ is bounded above, and by \inductionII, its denominator is bounded above by $r$, so the height of the fraction is bounded.  By \Height, therefore, 
$\dim V_{\lambda^2}$ is the value at $p_1$ of a Weyl polynomial which is a scalar multiple of $P^{R^1_1}_{\mu^1_1,\nu^1_1}(x)$.  By \Ratio, $\dim V_{\lambda^2_1} > \dim V_{\lambda^1_1}$, contrary to \inductionIII.

The induction step is similar, but we must exercise some additional care.  We write $G^1 = G^1_1\times H^1$.  We assume $p_1$ is large compared to the degree of the representation of $\dim V_{\lambda^1-\lambda^1_1}$.  
In other words, the upper bound on 
$$\eff(R^1_1)\log\dim V_{\lambda^1} - \lev(R^1_1)\log p_1$$
depends on $p_2,p_3,\ldots,p_{k+1}$, which we regard as constant, but not on $p_1$, which we take as large as we need.
Applying \Finiteness\ and \Height\ as before, there exists a simple factor $G^2_1$ of $G^2$, weights $\mu^2_1$ and $\nu^2_1$ of $G^2_1$, and a rational number $c\ge 1$, such that the root system $R^2_1$ has the same efficiency and level as $R^1_1$, $\lambda^2_1 = p\mu^2_1+\nu^2_1$, and
$$c P^{R^1_1}_{\mu^1_1,\nu^1_1}(x) = P^{R^2_1}_{\mu^2_1,\nu^2_1}(x).$$
Moreover $c=1$ implies $R^1_1 = R^2_1$.  As $\lambda^1_1$ and $\lambda^2_1$ are allowable, \Ratio\ guarantees that $\ord_{p_i}c = 0$ for all $i$.  Writing $G^2 = G^2_1\times H^2$, $\lambda^2 - \lambda^2_1$ is an allowable dominant weight of $H^2$, and for $i\ge 2$,
$$\eqalign{\ord_{p_i}\dim V_{\lambda^2-\lambda^2_1}
& = \ord_{p_i}\dim V_{\lambda^2}-\ord_{p_i}\dim V_{\lambda^2_1}\cr
& = \ord_{p_i}\dim V_{\lambda^1}-\ord_{p_i}P^{R_1^2}_{\mu_1^2,\nu_1^2}(p_1)\cr
& = \ord_{p_i}\dim V_{\lambda^1}-\ord_{p_i}P^{R_1^1}_{\mu_1^1,\nu_1^1}(p_1)\cr
& = \ord_{p_i}\dim V_{\lambda^2-\lambda^2_1} = \lev(R_i^1).\cr}$$
As $\dim V_{\lambda^2 - \lambda^2_1}$ divides ${\rm ht}(c)r\dim V_{\lambda^1-\lambda^1_1}$,
and ${\rm ht}(c)r$ can be bounded independent of the choice of primes $p_i$,
we can apply the induction hypothesis to conclude that
$$\dim V_{\lambda^2-\lambda^2_1} > \dim V_{\lambda^1-\lambda^1_1}.$$
Therefore
$$\eqalign{\dim V_{\lambda^2} 
&= \dim V_{\lambda^2_1}\dim V_{\lambda^2-\lambda^2_1} 
>  \dim V_{\lambda^2_1}\dim V_{\lambda^1-\lambda^1_1}\cr
&= P^{R^2_1}_{\mu^2_1,\nu^2_1}(p_1)\dim V_{\lambda^1-\lambda^1_1}
\ge P^{R^1_1}_{\mu^1_1,\nu^1_1}(p_1) \dim V_{\lambda^1-\lambda^1_1}\cr
&= \dim V_{\lambda^1_1} \dim V_{\lambda^1-\lambda^1_1}
= \dim V_{\lambda^1},\cr}$$
contrary to \inductionIII.
\qed

%---------------------------------------------------------------
\sec{The Gassmann Phenomenon}

Let $G=\SU(2)^n$ and $Z=Z(G)$.  Irreducible representations of $G$
are indexed by ordered $n$-tuples of non-negative integers: $\lambda=(a_1,\ldots,a_n)$.
The representation associated to $\lambda$ is
$$V_{\lambda}
=V_{(a_1,\ldots,0)}\otimes\cdots\otimes V_{(0,\ldots,a_n)}
=W_{a_1}\boxtimes\cdots\boxtimes W_{a_n},$$
where $W_i$ denotes the $i$th symmetric power of the natural representation of $\SU(2)$.
The restriction of $V_{\lambda}$ to $Z\cong\F_2^n$ is isotypic with character $\chi_{\lambda}$.
In coordinates,
$$\chi_{\lambda}(a_1,\ldots,a_n)=(-1)^{a_1z_1+\cdots+a_nz_n}.$$
For any character $\chi\in Z^*$, we write
$$X_\chi=\{\lambda\in\N^n\mid \chi_{\lambda}=\chi\}.$$
Thus, $X_\chi$ consists of $n$-tuples of non-negative integers
such that each coordinate has prescribed parity.

Let $Z_i$ ($i=1,2$) denote any subgroup of $Z$, and let $G^i = G/Z_i$.  Then,
$$\zeta_{G^i}(s)=\sum_{\chi\in Z_i^\perp}\sum_{\lambda\in X_\chi}(\dim V_{\lambda})^{-s}.$$
Now,
$$\dim V_{\lambda} = \prod_{i=1}^n (a_i+1),$$
so
$$\sum_{\lambda\in X_\chi}(\dim V_{\lambda})^{-s} = (2^{-s}\zeta(s))^O(\zeta(s)-2^{-s}\zeta(s))^E,$$
where $O$ (resp. $E$) is the number of coordinates which are constrained to be odd (resp. even)
in $X_\chi$.  Thus $O+E=n$ and $E-O$ is the trace of $\sigma(\chi)$, where
$\sigma\colon Z^*=\F_2^n\to\GL_n(\C)$ is defined by
$$\sigma(z_1^*,\ldots,z_n^*)=\left(\matrix{
(-1)^{z_1^*}&\cdots&0           \cr
\vdots      &\ddots&\vdots      \cr
0           &\cdots&(-1)^{z_n^*}\cr
}\right).$$
As long as
$$|\{z^*\in Z_1^\perp\mid\sigma(z^*)=c\}| = |\{z^*\in Z_2^\perp\mid\sigma(z^*)=c\}|$$
for all $c\in\Z$, we have $\zeta_{G^1}(s)=\zeta_{G^2}(s)$.

\nproc\Groups{Lemma}{There exists an integer $n$ and two injective homomorphisms 
$$\phi_1,\,\phi_2\colon\F_2^3\to\F_2^n$$
such that $\sigma(\phi_1(\F_2^3))$ and $\sigma(\phi_2(\F_2^3))$
are not conjugate in $\GL_n(\C)$ but for all $c\in\Z$,
$$|\{x\in\F_2^3\mid \tr(\sigma(\phi_1(x)))=c\}| 
= |\{x\in\F_2^3\mid \tr(\sigma(\phi_2(x)))=c\}|.$$
}

\proof
To guarantee non-conjugacy of images it suffices to arrange that the character of
$\sigma\scirc\phi_1$ precomposed by any automorphism of $\F_2^3$ is
distinct from the character of $\sigma\scirc\phi_2$.  To achieve this, we fix an injective
function 
$$f_1\colon \F_2^3\to 8\Z$$
such that
$$f_1(0)\ge \sum_{x\neq 0} |f(x)|.$$
By Fourier inversion, $f_1$ is the character of an effective representation
$\psi_1$.  Regarding it as a homomorphism to the diagonal matrices in $\GL_n(\C)$,
we deduce that there exists $\phi_1$ such that $\psi_1=\sigma\scirc\phi_1$.
Composing any non-linear automorphism of the pointed set $(\F_2^3,0)$ with $f_1$
we obtain a new function $f_2$, therefore a new representation $\psi_2$, and therefore a
new homomorphism $\phi_2$ satisfying the conditions of the lemma.
\qed

Setting $Z_i=\phi_i(\F_2^3)^\perp$ and $G^i=G/Z_i$, the lemma implies 
$\zeta_{G^1}(s)=\zeta_{G^2}(s)$.  Moreover, $G^1$ and $G^2$ cannot be isomorphic; indeed 
an isomorphism $\iota$ would give rise to an isomorphism $\tilde\iota$ between universal
covers, and therefore a diagram
$$\matrix{\llap{$G=$}\tilde G^1&\rt{\tilde\iota}&\tilde G^2\rlap{$=G$}\cr
\downarrow&&\downarrow\cr
G^1&\rt{\iota}&G^2\rlap.\cr}
$$
An inner automorphism is trivial on $Z$, so this means an element of $\Out(G)=S_n$
sends $Z_1$ to $Z_2$, contrary to \Groups.
We conclude:

\proc{Theorem}{There exist non-isomorphic compact semisimple Lie groups with
the same Witten zeta-function.}

\medskip

\biblio
\end